\documentclass[12pt]{amsart}

\newtheorem{theorem}{Theorem}[section]
\newtheorem{lemma}[theorem]{Lemma}
\newtheorem{corollary}[theorem]{Corollary}

\newtheorem{conjecture}[theorem]{Conjecture}

\theoremstyle{plain}
\newtheorem{definition}[theorem]{Definition}
\newtheorem{example}[theorem]{Example}
\newtheorem{remark}[theorem]{Remark}

\newtheorem{question}[theorem]{Question}

\newtheorem{setting}[theorem]{Setting}

\theoremstyle{definition}

\theoremstyle{remark}

\numberwithin{equation}{section}

\usepackage{fancyhdr}
\usepackage{amscd}
\usepackage{amsmath}
\usepackage{amsthm}
\usepackage{amssymb}

\usepackage[all]{xy}

\title[Pro-representability and Hodge numbers]{Pro-representability of Chow groups and Hodge numbers}

 \author{Sen Yang \\ \\ Dedicated to Kefeng Liu on his 61st birthday}

 \address{School of Mathematics and Finance \\  Chuzhou University \\
 Chuzhou, 239000, China\\
 }

\address{Applied Mathematics Research Center \\ Chuzhou University \\
 Chuzhou, 239000, China \\
}

 \email{yangsen.2014@tsinghua.org.cn, 101012424@seu.edu.cn}

\begin{document}

\date{}

\maketitle

\begin{abstract} 
Let $k$ be an algebraic field extension of $\mathbb{Q}$ and let $X$ be a smooth projective variety over $k$ of dimension $d \geq 2$. We study the pro-representability of the Chow group $CH^{p}(X)$ with $2 \leq p \leq d$. When certain Hodge numbers of $X$ vanish, namely, $H^{p}(X,\Omega^{i}_{X/k})=H^{p+1}(X,\Omega^{i}_{X/k})= \cdots =H^{2p-1-i}(X,\Omega^{i}_{X/k})=0$ for $i$ such that $0 \leq i \leq p-2$, we prove that the formal completion $\widehat{CH}^{p}(A)$ of $CH^{p}(X)$ at a local augmented Artinian $k$-algebra $A$ with the maximal ideal $m_{A}$ satisfies
\[
\widehat{CH}^{p}(A) \cong H^{p}(X, \Omega^{p-1}_{X/ k})\otimes_{k}m_{A}.
\]This provides a unified cohomological criterion for the pro-representability of the functor $\widehat{CH}^{p}$, generalizing earlier work by Bloch, Stienstra, and Mackall for $p=2$ and $p=3$. Our result reveals an intrinsic connection between the deformation theory of algebraic cycles and the Hodge structure of $X$.

\end{abstract}

\section{Introduction}
Let $X$ be a smooth projective variety over a field $k$ of characteristic zero. The Chow group $CH^{p}(X)$ of algebraic cycles of codimension-$p$ modulo rational equivalence is a fundamental invariant in algebraic geometry, intricately linked to the study of algebraic cycles, Hodge theory, and arithmetic properties of varieties. A central and subtle problem is understanding how the Chow groups behave under infinitesimal deformations of the underlying variety $X$. This paper is devoted to the study of the pro-representability of Chow groups.

Let $Art_{k}$ denote the category of augmented local Artinian $k$-algebras. For $A \in Art_{k}$, we write $X_A=X\times_{\mathrm{Spec}(k)}\mathrm{Spec}(A)$. Let $K^{M}_{p}(O_{X_{A}})$ be the Milnor K-theory sheaf associated to the presheaf
\[
 U \to K^{M}_{p}(O_{X_{A}}(U)),
\] where $U = \mathrm{Spec}(R) \subset X$ is open affine, $K^{M}_{p}(O_{X_{A}}(U))=K^{M}_{p}(R \otimes_{k}A)$. When $A=k$, we write $K^{M}_{p}(O_{X})=K^{M}_{p}(O_{X_{k}})$.

A fundamental tool in this study is the Bloch formula (cf. \cite{Bl2-Annals, Ka,Kerz,Quillen,Soule}). Since we are interested in Chow groups with $\mathbb{Q}$-coefficients, we use the following format of Bloch formula
\begin{equation} \label{eq:Bloch formula}
CH^{p}(X)\otimes_{\mathbb{Z}}\mathbb{Q} \cong H^{p}(X, K^{M}_{p}(O_{X}))\otimes_{\mathbb{Z}}\mathbb{Q}.
\end{equation}

Bloch formula \eqref{eq:Bloch formula} motivates a functor on $Art_{k}$
{\footnotesize
\begin{align} \label{eq:completion}
\widehat{CH}^{p}: \ A \to \mathrm{ker} \ \{ H^{p}(X, K^{M}_{p}(O_{X_A}))\otimes_{\mathbb{Z}}\mathbb{Q}  \xrightarrow{aug} H^{p}(X, K^{M}_{p}(O_{X}))\otimes_{\mathbb{Z}}\mathbb{Q} \},
\end{align}
}where $H^{p}(X, K^{M}_{p}(O_{X_A}))\otimes_{\mathbb{Z}}\mathbb{Q}$ may be considered as an infinitesimal deformation of $CH^{p}(X)\otimes_{\mathbb{Z}}\mathbb{Q}$ via the isomorphism \eqref{eq:Bloch formula} and the map $aug$ is induced by the augmentation $A \to k$. 

The group $\widehat{CH}^{p}(A)$ is called the formal completion of $CH^{p}(X)$. It is clear that $\widehat{CH}^{p}$ is a functor of Artin rings. A functor of Artin rings $F$ is called pro-representable if it is isomorphic to the functor $h_{R}$, which is defined as 
\[
h_{R}(A)=Hom(R,A).
\]Here $A \in  Art_{k}$, $R$ is a complete local Noetherian $k$-algebra with the maximal ideal $m_{R}$ such that $R/m_{R}^{n} \in  Art_{k}$, and $Hom(R,A)$ denotes the set of local $k$-algebra homomorphisms from $R$ to $A$.

A useful example of pro-representable functor is given as
\begin{equation} \label{eq: pro-rep vs}
A \to V \otimes_{k}m_{A},
\end{equation}where $V$ is a finite dimensional $k$-vector space and $m_{A}$ is the maximal ideal of $A \in Art_{k}$. We refer the readers to \cite{Sc} for details on the theory of pro-representable functors.

We are interested in the pro-representability of the functor $\widehat{CH}^{p}$. When $p=1$, it is well-known that 
\begin{equation} \label{eq:ch1}
\widehat{CH}^{1}(A) \cong H^{1}(X,O_{X})\otimes_{k} m_{A},
\end{equation}which carries the form of \eqref{eq: pro-rep vs}. Hence, the pro-representability of $\widehat{CH}^{1}$ is clear.

The study of the pro-representability of $\widehat{CH}^{p}(A)$ is highly non-trivial for $p \geq 2$. In the case of $p=2$, when $X$ is a surface defined over a field $k$, which is an algebraic field extension of $\mathbb{Q}$, Bloch \cite{Bl3, Bl4} proves that there is a short exact sequence
\[
0 \to H^{2}(X,O_{X})\otimes_{k} \frac{\Omega^{1}_{A/k}}{dm_{A}} \to \widehat{CH}^{2}(A) \to H^{2}(X,\Omega^{1}_{X/k})\otimes_{k} m_{A} \to 0.
\]This yields the following result.
\begin{theorem}[\cite{Bl3}] \label{t:Bloch}
Let $X$ be a smooth projective surface over an algebraic field extension $k$ of $\mathbb{Q}$, and suppose that $ H^{2}(X,O_{X})=0$. Then, for any $A \in Art_{k}$, there is an isomorphism
\[
\widehat{CH}^{2}(A) \cong H^{2}(X,\Omega^{1}_{X/k})\otimes_{k} m_{A}.
\]Consequently, the functor $\widehat{CH}^{2}$ is pro-representable, cf. \eqref{eq: pro-rep vs}.

\end{theorem}

This motivated the following important conjecture, which has been intensively studied; see \cite{BKL,BS,Hu,PW,V2,V3} for progress on it. Recall that the geometric genus $p_{g}(X)=\mathrm{dim} H^{0}(X,\Omega^{2}_{X/k})=\mathrm{dim} H^{2}(X,O_{X})$.
\begin{conjecture} [\cite{Bl3}] \label{c:Bloch}
Let $X$ be a connected smooth complex projective surface with trivial geometric genus, i.e., $p_{g}(X)=0$. Then the Albanese map 
\[
CH^{2}_{deg \ 0}(X) \to Alb(X)
\]
is an isomorphism, where $CH^{2}_{deg \ 0}(X)$ is the subgroup of $CH^{2}(X)$ consisting of zero cycles with degree zero and $Alb(X)$ is the Albanese variety.
\end{conjecture}

Very recently, Mackall extended Theorem \ref{t:Bloch} to the case of $p=3$.
\begin{theorem} [\cite{Ma}] \label{t:mackall}
Let $X$ be a smooth, projective, geometrically connectd variety over an algebraic field extension $k$ of $\mathbb{Q}$. Suppose that the following conditions are satisfied:
\begin{align*}\label{eq: condi-vanish-explicit}
\begin{cases}
   H^{3}(X,O_{X})= H^{4}(X,O_{X})= H^{5}(X,O_{X})=0,  \\
  H^{3}(X,\Omega^{1}_{X/k})= H^{4}(X, \Omega^{1}_{X/k})=0.
\end{cases}
\end{align*}Then, for any $A \in Art_{k}$, there is an isomorphism
\[
\widehat{CH}^{3}(A) \cong H^{3}(X,\Omega^{2}_{X/k})\otimes_{k} m_{A}.
\]Consequently, the functor $\widehat{CH}^{3}$ is pro-representable,  cf. \eqref{eq: pro-rep vs}.

\end{theorem}

For cases $p > 3$, what can be said about the pro-representability of the functor $\widehat{CH}^{p}$? This leads naturally to the following question.
\begin{question} \label{q:main ques}
Let $X$ be a smooth projective variety over a field $k$ of dimension $d \geq 2$, where $k$ is an algebraic field extension of $\mathbb{Q}$. For $p$ an integer satisfying that $2 \leq p \leq d$, when is the functor $\widehat{CH}^{p}$ pro-representable ?

\end{question}

The work of Bloch and Mackall suggests a profound and general principle: the pro-representability of $\widehat{CH}^{p}$ is governed by vanishing of certain Hodge numbers of $X$. Following this idea, we provide an answer to Question \ref{q:main ques}, which is the principal result of this paper.
\begin{theorem} [cf. Theorem \ref{t:yang} below]\label{t:yang1}
Let $X$ be a smooth projective variety over a field $k$ of dimensin $d \geq 2$, where $k$ is an algebraic field extension of $\mathbb{Q}$. Suppose that the following conditions are satisfied:
\begin{equation*} \label{eq: condi-vanish}
H^{p}(X,\Omega^{i}_{X/k})= H^{p+1}(X,\Omega^{i}_{X/k})= \cdots =H^{2p-1-i}(X,\Omega^{i}_{X/k})=0,
\end{equation*}where $2 \leq p \leq d$ and $i$ is an integer satisfying that $0 \leq i \leq p-2$.
Then, for any $A \in Art_{k}$, there is an isomorphism
\[
\widehat{CH}^{p}(A) \cong H^{p}(X, \Omega^{p-1}_{X/ k})\otimes_{k}m_{A}.
\]Consequently, the functor $\widehat{CH}^{p}$ is pro-representable,  cf. \eqref{eq: pro-rep vs}.

\end{theorem}

It is worth noting that in Theorem \ref{t:Bloch}, Theorem \ref{t:mackall} and Theorem \ref{t:yang1}, the ground field $k$ is assumed to be an algebraic field extension of $\mathbb{Q}$. This is to circumvent difficulties in computing $\widehat{CH}^{p}(A)$ that arise when the differential $\Omega^{1}_{k/ \mathbb{Q}} \neq 0$. 

If the ground field $k$ contains transcendental elements, the computation of $\widehat{CH}^{p}(A)$ (even for $p=2$) is much more complicated because of the appearance of non-trivial absolute differentials. When $p=2$, this is studied by Stienstra \cite{Stien1}.
  
This paper is structured as follows. In Section 2 and Section 3, we recall the necessary background on de Rham cohomology and cyclic homology. We prove our main result, Theorem \ref{t:yang} in Section 4.

\textbf{Notations:} 

(1) If not stated otherwise, $k$ is a field of characteristic zero. Let $Art_{k}$ denote the category of local Artinian $k$-algebra $A$ with the maximal ideal $m_{A}$ such that $A/m_{A}=k$. 


(2) For any abelian group $M$, $M_{\mathbb{Q}}$ denotes $M \otimes_{\mathbb{Z}} \mathbb{Q}$.

\section{de Rham cohomology}

For $X$ a scheme of finite type over $k$, we denote by $\Omega^{\bullet}_{X/k}$ the de Rham complex
\[
O_{X} \xrightarrow{d} \Omega^{1}_{X/k} \xrightarrow{d} \Omega^{2}_{X/k} \xrightarrow{d} \cdots.
\]For each non-negative integer $j$, the de Rham cohomology $H^{j}_{dR}(X/k)$ is defined to be the hypercohomology of the de Rham complex $\Omega^{\bullet}_{X/k}$, i.e.,
\[
H^{j}_{dR}(X/k)=\mathbb{H}^{j}(X, \Omega^{\bullet}_{X/k}).
\]

We recall a concrete example of de Rham cohomology.
\begin{example} [cf. remark from page 54 of \cite{Ha}]
The de Rham cohomology of the affine line $A^{1}_{k}$ is as following:
\begin{align*} \label{eq:de rhamline}
\begin{cases}
 \begin{CD}
   & H^{j}_{dR}(A^{1}_{k}/k)=k, \  j=0; \\
  &H^{j}_{dR}(A^{1}_{k}/k)=0, \ j \neq 0.
 \end{CD}
\end{cases}
\end{align*}

\end{example}

\begin{lemma}[\cite{Gro1}]  \label{l:affineDerham}
When $X=\mathrm{Spec}(T)$ is affine, the de Rham cohomology $H^{j}_{dR}(X/k)$ agrees with the $j$-th cohomology of the following cochain complex:
\[
\Omega^{\bullet}_{T/k}: \ T \xrightarrow{d} \Omega^{1}_{T /k} \xrightarrow{d} \Omega^{2}_{T /k} \xrightarrow{d} \cdots \xrightarrow{d} \Omega^{j}_{T /k} \xrightarrow{d} \cdots.
\]

\end{lemma}

We sketch the proof briefly. 
\begin{proof}
There is a spectral sequence
\[
   E_{1}^{i_{1},i_{2}}=H^{i_{2}}(X, \Omega^{i_{1}}_{X/k}) \Longrightarrow {H}_{dR}^{j}(X/k).
\]where $j=i_{1}+i_{2}$. By Serre vanishing theorem, when $i_{2}>0$, for each $\Omega^{i_{1}}_{X/k}$, we have
\[
H^{i_{2}}(X, \Omega^{i_{1}}_{X/k})=0.
\]This implies that ${H}_{dR}^{j}(X/k)$ is the $j$-th cohomology of the following complex 
{\footnotesize
\[
H^{0}(X, O_{X}) \to  H^{0}(X, \Omega^{1}_{X/k}) \to \cdots \to H^{0}(X, \Omega^{j}_{X/k}) \to \cdots,
\]
}which is of the form
\[
T \xrightarrow{d} \Omega^{1}_{T /k} \xrightarrow{d} \Omega^{2}_{T /k} \xrightarrow{d} \cdots \xrightarrow{d} \Omega^{j}_{T /k} \xrightarrow{d} \cdots.
\]
\end{proof}

De Rham cohomology enjoys a list of nice properties. We refer to \cite{Gro1,Gro2,Ha,The} for details.

Let $\mathcal{H}_{dR}^{q}(O_{X})$ be the de Rham cohomology sheaf, which is defined to be the sheaf associated to the presheaf
\[
U \to H_{dR}^{q}(U/k),
\]where $U=\mathrm{Spec}(R) \subset X$ is open affine, and ${H}_{dR}^{q}(U/k)$ is the hypercohomology of the de Rham complex
\[
\Omega^{\bullet}_{U/k}: \ U \xrightarrow{d} \Omega^{1}_{U /k} \xrightarrow{d} \Omega^{2}_{U /k} \xrightarrow{d} \cdots.
\]

We are interested in finding a resolution of the sheaf $\mathcal{H}_{dR}^{q}(O_{X})$, which enables us to compute its cohomology. This leads to the classical Bloch-Ogus theorem, originally proved by Bloch and Ogus \cite{BO}.

Given a smooth algebraic variety $X$ and a cohomology theory $h$ satisfying natural axioms, the classical Bloch-Ogus theorem, says that the Zariski sheafification of the Cousin complex (formed from the coniveau spectral sequence) of $h$ is a flasque resolution of the Zariski sheaf associated to the presheaf $U \to h^{\ast}(U)$.

Bloch and Ogus \cite{BO} proved their theorem for \'etale cohomology with coefficients in roots of unity, by reducing to the effacement theorem. Gabber extended it in \cite{Gabber}. Colliot-Th\'el\`ene, Hoobler and Kahn applied Gabber's method to prove the theorem for any cohomology theory with support which satisfies \'etale excision and homotopy invariance (or projective bundle formula), see Corollary 5.1.11 and Proposition 5.4.3 of \cite{CTHK}.
We have applied the Bloch-Ogus theorem to algebraic K-theory and cyclic homology to study deformations of Chow groups in \cite{DHY,Y5}. We will apply it to the de Rham cohomology (see Example 2.2 of \cite{BO} or 7.3 of \cite{CTHK}) in this paper.

In the remainder of this section, $X$ is a smooth projective variety over a field $k$ of characteristic zero, and $q$ is a non-negative integer.
\begin{theorem} [cf. Theorem 4.2 of \cite{BO}]  \label{t:BGQresolution}
The Zariski sheafification of the Cousin complex
{\footnotesize
\begin{align} \label{eq: BGQ seq1}
0  \rightarrow  \bigoplus_{x \in X^{(0)}}H^{q}_{dR}(x/k) \to \bigoplus_{x \in X^{(1)}}H^{q-1}_{dR}(x/k) \to \bigoplus_{x \in X^{(2)}}H^{q-2}_{dR}(x/k) \to \cdots,
\end{align}
}is a flasque resolution of the de Rham cohomology sheaf $\mathcal{H}_{dR}^{q}(O_{X})$, where each $H^{m}_{dR}(x/k)$, $m=q, q-1, \cdots$, is defined as the following
\[
H^{m}_{dR}(x/k)=\varinjlim_{x \in U}H^{m}_{dR}(U/k).
\]

\end{theorem}

\begin{corollary} [cf. Corollary 6.2 of \cite{BO}]  \label{c:BO}
For each integer $j$ such that $j > q \geq 0$, the $j$-th cohomology of the de Rham cohomology sheaf $\mathcal{H}_{dR}^{q}(O_{X})$ vanishes, i.e.,
\[
H^{j}(X,\mathcal{H}_{dR}^{q}(O_{X}))=0.
\]
\end{corollary}

\begin{proof}
By Theorem \ref{t:BGQresolution}, the sheaf cohomology $H^{j}(X,\mathcal{H}_{dR}^{q}(O_{X}))$ can be computed as the $j$-th cohomology of the complex \eqref{eq: BGQ seq1}. To be precise, it is the $j$-th cohomology of the complex
{\footnotesize
\begin{align*} \label{eq: BGQ seq2}
& 0  \rightarrow  \bigoplus_{x \in X^{(0)}}H^{q}_{dR}(x/k) \to \cdots \to \bigoplus_{x \in X^{(j-1)}}H^{q-(j-1)}_{dR}(x/k) \to  \bigoplus_{x \in X^{(j)}}H^{q-j}_{dR}(x/k) \\
& \to \bigoplus_{x \in X^{(j+1)}}H^{q-(j+1)}_{dR}(x/k) \to \cdots.  \notag
\end{align*}
}

  Since  $j > q$, $q-j < 0$, $H^{q-j}_{dR}(x/k)=0$. This yields the desired result.

\end{proof}

For $A \in Art_{k}$, we write $X_{A}=X \otimes_{\mathrm{Spec}(k)} \mathrm{Spec}(A)$. Let $\mathcal{H}_{dR}^{q}(O_{X_{A}})$ be the de Rham cohomology sheaf, which is defined to be the sheaf associated to the presheaf
\[
U \to H_{dR}^{q}(U_{A}/k),
\]where $U=\mathrm{Spec}(R) \subset X$ is open affine, $U_{A}=U \times_{\mathrm{Spec}(k)} \mathrm{Spec}(A)$ and ${H}_{dR}^{q}(U_{A}/k)$ is the hypercohomology of the de Rham complex $\Omega^{\bullet}_{U_{A}/k}$, which has the form:
\[
U_{A}\xrightarrow{d} \Omega^{1}_{U_{A}/k} \xrightarrow{d} \Omega^{2}_{U_{A}/k} \xrightarrow{d} \cdots.
\]

\begin{lemma}[K$\ddot{\mathrm{u}}$nneth formula] \label{l:1}
With the notation as above, there is an isomorphism between sheaves
\begin{equation} \label{eq:ku1}
\mathcal{H}_{dR}^{q}(O_{X_{A}}) \cong \bigoplus_{q_{1}+q_{2}=q}\mathcal{H}_{dR}^{q_{1}}(O_{X})\otimes_{k} H^{q_{2}}_{dR}(A),
\end{equation}where $q$ is a non-negative integer, and $H^{q_{2}}_{dR}(A)$ is the de Rham cohomology of $A$, which is the $q_{2}$-th cohomology of the following complex (cf. Lemma \ref{l:affineDerham})
\[
\Omega^{\bullet}_{A/k}:\ A \xrightarrow{d} \Omega^{1}_{A/k} \xrightarrow{d} \Omega^{2}_{A/k} \xrightarrow{d} \cdots.
\]

\end{lemma}

The K$\ddot{\mathrm{u}}$nneth formula has been known to experts. We sketch a proof for the reader's convenience and refer to section 8 of \cite{The} for a comprehensive study of K$\ddot{\mathrm{u}}$nneth formulas of de Rham cohomology.

\begin{proof}
It suffices to check the isomorphism \eqref{eq:ku1} locally on $U=\mathrm{Spec}(R) \subset X$ and show the following isomorphism:
\[
H_{dR}^{q}(U_{A}/k) \cong \bigoplus_{q_{1}+q_{2}=q} H_{dR}^{q_{1}}(U/k)\otimes_{k} H^{q_{2}}_{dR}(A).
\]

Since both $U$ and $U_{A}=\mathrm{Spec}(R \otimes_{k} A)$ are affine, by Lemma \ref{l:affineDerham}, the de Rham cohomology $H_{dR}^{q}(U_{A}/k)$ and $H_{dR}^{q_{1}}(U/k)$ agree with the cohomology groups of the following complexes $\Omega^{\bullet}_{R\otimes_{k}A /k}$ and $\Omega^{\bullet}_{R/k}$ respectively,
\begin{align*}
 \Omega^{\bullet}_{R \otimes_{k}A /k}&:  R \otimes_{k}A \xrightarrow{d} \Omega^{1}_{R \otimes_{k}A /k} \xrightarrow{d} \Omega^{2}_{R \otimes_{k}A /k} \xrightarrow{d} \cdots, \\
\Omega^{\bullet}_{R/k}& :  R \xrightarrow{d} \Omega^{1}_{R/k} \xrightarrow{d} \Omega^{2}_{R/k} \xrightarrow{d} \cdots.
\end{align*}
That is, we have 
\[
H_{dR}^{q}(U_{A}/k) \cong H^{q}(\Omega^{\bullet}_{R \otimes_{k}A /k}), \  \ H^{q_{1}}_{dR}(U /k) \cong H^{q_{1}}(\Omega^{\bullet}_{R /k}).
\]

As a consequence, it is sufficient to show that 
\begin{equation*} 
H^{q}(\Omega^{\bullet}_{R \otimes_{k} A/k}) \cong \bigoplus_{q_{1}+q_{2}=q}H^{q_{1}}(\Omega^{\bullet}_{R/k})\otimes_{k} H^{q_{2}}(\Omega^{\bullet}_{A/k}).
\end{equation*}

To achieve this, we note that the following isomorphism
\[
\Omega^{1}_{R \otimes_{k} A/k} \cong (\Omega^{1}_{R/k} \otimes_{k}A) \oplus (\Omega^{1}_{A/k} \otimes_{k}R),
\]implies that the complex $\Omega^{\bullet}_{R \otimes_{k} A/k}$ can be identified with the complex associated to the tensor product $\Omega^{\bullet}_{R/k} \otimes_{k} \Omega^{\bullet}_{A/k}$. This yields an isomorphism between cohomology groups of complexes:
\[
H^{q}(\Omega^{\bullet}_{R \otimes_{k} A/k}) \cong  H^{q}(\Omega^{\bullet}_{R/k} \otimes_{k} \Omega^{\bullet}_{A/k}).
\]

It remains to shows that 
\begin{equation} \label{eq:1}
H^{q}(\Omega^{\bullet}_{R/k} \otimes_{k} \Omega^{\bullet}_{A/k}) \cong \bigoplus_{q_{1}+q_{2}=q}H^{q_{1}}(\Omega^{\bullet}_{R/k})\otimes_{k} H^{q_{2}}(\Omega^{\bullet}_{A/k}).
\end{equation}

There is a K$\ddot{\mathrm{u}}$nneth spectral sequence (cf. Theorem 2.20 of \cite{Mc})
\[
E^{p,q}_{2}=\bigoplus_{q_{1}+q_{2}=q}Tor^{p}_{k}(H^{q_{1}}(\Omega^{\bullet}_{R/k}), H^{q_{2}}(\Omega^{\bullet}_{A/k})) \Longrightarrow H^{p+q}(\Omega^{\bullet}_{R \otimes_{k} A/k}).
\]

Since $k$ is a field, there is no non-trivial higher Tor groups. To be precise, when $p>0$, 
\[
Tor^{p}_{k}(H^{q_{1}}(\Omega^{\bullet}_{R/k}), H^{q_{2}}(\Omega^{\bullet}_{A/k}))=0.
\]Hence, the above spectral sequence degenerates and yields the isomorphism \eqref{eq:1}. This completes the proof.

\end{proof}

Let $j$ be a non-negative integer. We pass the isomorphism \eqref{eq:ku1} to cohomology groups:
\begin{align*}
H^{j}(X,\mathcal{H}_{dR}^{q}(O_{X_{A}})) & \cong H^{j}(X, \bigoplus_{q_{1}+q_{2}=q} \mathcal{H}_{dR}^{q_{1}}(O_{X})\otimes_{k} H_{dR}^{q_{2}}(A)) \\
& \cong  \bigoplus_{q_{1}+q_{2}=q} H^{j}(X, \mathcal{H}_{dR}^{q_{1}}(O_{X})\otimes_{k} H_{dR}^{q_{2}}(A)) \\
& \cong  \bigoplus_{q_{1}+q_{2}=q} H^{j}(X, \mathcal{H}_{dR}^{q_{1}}(O_{X}))\otimes_{k} H_{dR}^{q_{2}}(A).
\end{align*}

When $j>q$, since $q \geq q_{1}$, it is obvious that $j \geq q_{1}$. It is immediately from Corollary \ref{c:BO} that $ H^{j}(X, \mathcal{H}_{dR}^{q_{1}}(O_{X}))=0$. This shows the following result.
\begin{corollary} \label{c:BO-A}
For each integer $j$ such that $j > q \geq 0$, we have
\[
H^{j}(X,\mathcal{H}_{dR}^{q}(O_{X_{A}}))=0.
\]

\end{corollary}

By Lemma \ref{l:affineDerham}, the de Rham cohomology sheaf $\mathcal{H}_{dR}^{q}(O_{X_{A}})$ can be identified with the $q$-th cohomology of the de Rham complex
\[
\cdots \xrightarrow{d} \Omega^{q-1}_{X_{A}/k} \xrightarrow{d} \Omega^{q}_{X_{A}/k} \xrightarrow{d}  \Omega^{q+1}_{X_{A}/k} \xrightarrow{d} \cdots.
\]To be precise, let $T^{q}$ denote the kernel of 
\[
d: \ \Omega^{q}_{X_{A}/k} \to \Omega^{q+1}_{X_{A}/k},
\]and let $d\Omega^{q-1}_{X_{A}/k}$ denote the image of 
\[
d: \ \Omega^{q-1}_{X_{A}/k} \to \Omega^{q}_{X_{A}/k}.
\]Then the de Rham cohomology sheaf $\mathcal{H}_{dR}^{q}(O_{X_{A}})$ can be identified with the quotient sheaf $\frac{T^{q}}{d\Omega^{q-1}_{X_{A}/k}}$, which should be known to experts. 
\begin{lemma} \label{l:derhshiso}
With the notation as above, there is an isomorphism
\[
\mathcal{H}_{dR}^{q}(O_{X_{A}}) \cong \frac{T^{q}}{d\Omega^{q-1}_{X_{A}/k}}.
\]
\end{lemma}

We use Corollary \ref{c:BO-A} and Lemma \ref{l:derhshiso} to deduce the following result.
\begin{corollary} \label{c:BGQresolution}
For each integer $j$ such that $j > q \geq 0$, the $j$-th cohomology of the sheaf $\frac{T^{q}}{d\Omega^{q-1}_{X_{A}/k}}$ vanishes, i.e.,
\[
H^{j}(X, \frac{T^{q}}{d\Omega^{q-1}_{X_{A}/k}})=0.
\]
\end{corollary}

\section{Cyclic homology}

Let $R$ be an associative commutative unital $k$-algebra. The Hochschild complex (with coefficients in $R$) is 
\begin{equation*}
C(R): \ \ \cdots \to R^{\otimes (n+1)} \xrightarrow{b} R^{\otimes n} \to \cdots \to R\otimes_{k}R \xrightarrow{b} R,
\end{equation*} in which $b$ is the Hochschild boundary (cf. section 1.1.1 of \cite{Loday}). The Hochschild homology group $HH_{n}(R)$ is defined to be the homology of the Hochschild complex $C(R)$
\[
HH_{n}(R):=H_{n}(C(R)).
\]

The cyclic bicomplex $CC(R)$ is a first-quadrant bicomplex formed from the Hochschild complex $C(R)$, whose total complex $\mathrm{Tot}(CC(R))$ can be used to define the cyclic homology of $R$
\begin{equation} \label{eq:hc def}
HC_{n}(R):=H_{n}(\mathrm{Tot}(CC(R))).
\end{equation} We refer to section 5.1 of \cite{Loday} for details.

Weibel, McCarthy and Keller generalized the definition of cyclic homology and its variants to schemes, exact categories and localization pairs; see section 2 of \cite{CHSW} for a summary.

Lambda operations $\lambda^{m}$ and Adams operations $\psi^{m}$ are defined on Hochschild homology $HH_{n}(R)$ and cyclic homology $HC_{n}(R)$, see section 4.5 of \cite{Loday} and section 9.4.3 of \cite{W2} for details. For each integer $n \geq 1$, the group $HC_{n}(R)$ decomposes into direct sums of eigenspaces
\begin{equation*}
HC_{n}(R)=HC^{(1)}_{n}(R) \oplus \cdots \oplus HC^{(i)}_{n}(R) \oplus \cdots \oplus HC^{(n)}_{n}(R),
\end{equation*}where $HC^{(i)}_{n}(R)$ is the eigenspace of $\psi^{m}=m^{i+1}$. For $n=0$, $HC_{0}(R)=HC^{(0)}_{0}(R)=R$.

\begin{lemma} [cf. Corollary 9.8.16 of \cite{W2}]\label{lemma: l-l-omega}
With the notation as above, there is an isomorphism
\[
HC^{(n)}_{n}(R) \cong \dfrac{\Omega^{n}_{R/k}}{d\Omega^{n-1}_{R/k}},
\]where $d\Omega^{n-1}_{R/k}$ is the image of $d: \Omega^{n-1}_{R/k} \to \Omega^{n}_{R/k}$.
\end{lemma}

For $A \in Art_{k}$, the augmentation $A \to k$ induces a morphism of algebraic K-groups
\[
K_{n}(R\otimes_{k}A) \to  K_{n}(R\otimes_{k}k) \xrightarrow{\cong } K_{n}(R),
\]whose kernel is the relative K-group, denoted by $\overline{K}_{n}(R\otimes_{k}A)$. By \cite{Soule}, these K-groups carry Adams operations $\psi^{m}$. Let $\overline{K}_{n}^{(i)}(R\otimes_{k}A)$ be the eigenspace of $\psi^{m}=m^{i}$. The relative cyclic homology $\overline{HC}_{n-1}(R\otimes_{k}A)$ and its eigenspace $\overline{HC}_{n-1}^{(i-1)}(R\otimes_{k}A)$ are defined similarly. Goodwillie and Cathelineau proved that these relative groups are connected by the relative Chern character.
\begin{theorem} [\cite{Cath,Good}] \label{theorem: GoodwillieCathelineau}
With the notation as above, the relative Chern character induces the following isomorphisms $\mathrm{(}$ where cyclic homology is defined over $\mathbb{Q} \ \mathrm{)}$,
 \begin{align*}\label{eq: relchern}
\begin{cases}
   \overline{K}_{n}(R\otimes_{k}A)_{\mathbb{Q}} \xrightarrow{\cong} \overline{HC}_{n-1}(R\otimes_{k}A),\\
  \overline{K}_{n}^{(i)}(R\otimes_{k}A)_{\mathbb{Q}}  \xrightarrow{\cong} \overline{HC}_{n-1}^{(i-1)}(R\otimes_{k}A).
\end{cases}
\end{align*}

\end{theorem}

This theorem is very useful for computing relative K-groups. Corti$\mathrm{\tilde{n}}$as-Haesemeyer-Weibel \cite{CHW} generalized it to the space level.

In the remainder of this section, $X$ is a smooth projective variety over a field $k$ of characteristic zero. For $A \in Art_{k}$, we write $X_{A}=X \times \mathrm{Spec}(A)$.
\begin{definition}

Let $HC_{n}(O_{X_{A}})$ be the cyclic homology sheaf associated to the following presheaf
\[
 U \to HC_{n}(O_{X_{A}}(U)), 
\]
where $U=\mathrm{Spec}(R) \subset X$ is open affine, $HC_{n}(O_{X_{A}}(U))=HC_{n}(R\otimes_{k}A)$ is the cyclic homology of $R\otimes_{k}A$, cf \eqref{eq:hc def}. When $A=k$, we write $HC_{n}(O_{X})=HC_{n}(O_{X_{k}})$.

\end{definition}

Adams operations $\psi^{m}$ act on the sheaf $HC_{n}(O_{X_{A}})$, whose eigenspaces (of $\psi^{m}=m^{i+1}$) are $HC^{(i)}_{n}(O_{X_{A}})$. When $A=k$, we write $HC^{(i)}_{n}(O_{X})= HC^{(i)}_{n}(O_{X_{k}})$.

\begin{definition} \label{d:relcyc}
 The relative cyclic homology sheaf $\overline{HC}^{(i)}_{n}(O_{X_{A}})$ is defined to be the kernel of the following map (induced by $A \to k$) 
\[
 HC^{(i)}_{n}(O_{X_{A}}) \to HC^{(i)}_{n}(O_{X}).
\]

\end{definition}

Let $p$ be an integer such that $1 \leq p \leq \mathrm{dim}(X)$. We denote by $\overline{K}_{p}^{M}(O_{X_{A}})$ the relative Milnor K-theory sheaf, which is defined as the kernel of the following map
\[
K_{p}^{M}(O_{X_{A}}) \xrightarrow{A \to k} K_{p}^{M}(O_{X}).
\]

There is a short exact sequence 
\[
0 \to \overline{K}_{p}^{M}(O_{X_{A}}) \to K_{p}^{M}(O_{X_{A}}) \to K_{p}^{M}(O_{X}) \to 0,
\]which is split. It follows that the formal completion $\widehat{CH}^{p}(A)$, which is defined in \eqref{eq:completion}, is identified with $H^{p}(X,\overline{K}_{p}^{M}(O_{X_{A}}))_{\mathbb{Q}}$, i.e.,
\begin{equation} \label{eq:relK}
\widehat{CH}^{p}(A)  \cong H^{p}(X,\overline{K}_{p}^{M}(O_{X_{A}}))_{\mathbb{Q}}.
\end{equation}

Let $K_{p}^{(p)}(O_{X_{A}})$ be the eigenspace of $\psi^{m}=m^{p}$, where $\psi^{m}$ is Adams operation on the K-theory sheaf $K_{p}(O_{X_{A}})$, cf. \cite{Soule}. By a sheaf version of Theorem 5 in page 526 of \cite{Soule}, there exists an isomorphism $K_{p}^{M}(O_{X_{A}})_{\mathbb{Q}} \xrightarrow{\cong} K_{p}^{(p)}(O_{X_{A}})_{\mathbb{Q}}$, which yields that 
\begin{equation} \label{eq:releigen}
\overline{K}_{p}^{M}(O_{X_{A}})_{\mathbb{Q}}  \xrightarrow{\cong} \overline{K}_{p}^{(p)}(O_{X_{A}})_{\mathbb{Q}},
\end{equation}where $\overline{K}_{p}^{(p)}(O_{X_{A}})$ is the kernel of the map $K_{p}^{(p)}(O_{X_{A}}) \xrightarrow{A \to k} K_{p}^{(p)}(O_{X})$.

By combining a sheaf version of the isomorphism of Theorem \ref{theorem: GoodwillieCathelineau} (let $i=n=p$) with \eqref{eq:relK} and \eqref{eq:releigen}, we can compute the formal completion $\widehat{CH}^{p}(A)$ in terms of cyclic homology.

\begin{lemma} \label{l:chhatisocyc}
With the notation as above, there is an isomorphism
\begin{equation*}
\widehat{CH}^{p}(A) \cong H^{p}(X,\overline{HC}^{(p-1)}_{p-1}(O_{X_{A}})).
\end{equation*}

\end{lemma}

\section{Pro-representability of Chow groups}

In this section, we study the pro-representability of $\widehat{CH}^{p}(A)$ defined in \eqref{eq:completion}, where $1 \leq p \leq \mathrm{dim}(X)$. Since the case of $p=1$ is known, cf. \eqref{eq:ch1}, we focus on the cases of $p \geq 2$. This requires that $\mathrm{dim}(X) \geq 2$. 

From now on, we adopt the following setting.

\begin{setting} \label{s:set1}
Let $X$ be a smooth projective variety over a field $k$ of characteristic zero, which has dimension $\mathrm{dim}(X) \geq 2$. Let $p$ be an integer such that $2 \leq p \leq \mathrm{dim}(X)$. We fix the integer $p$ and assume that $X$ satisfies the following condition
\begin{equation} \label{eq: condi-vanish}
H^{p}(X,\Omega^{i}_{X/k})= H^{p+1}(X,\Omega^{i}_{X/k})= \cdots =H^{2p-1-i}(X,\Omega^{i}_{X/k})=0,
\end{equation}where $i$ is an integer satisfying that $0 \leq i \leq p-2$.

\end{setting}

Recall that the Hodge number $h^{i,j}(X)$ is the dimension of $H^{j}(X,\Omega^{i}_{X/k})$, where $j$ is a non-negative integer. The condition \eqref{eq: condi-vanish} can be alternatively stated as:
\[
h^{i,p}(X)=h^{i,p+1}(X)= \cdots = h^{i,2p-1-i}(X)=0.
\]

To understand the condition \eqref{eq: condi-vanish}, we first look at the case of $p=2$. Note that $i$ can only be $0$ in this case. The condition \eqref{eq: condi-vanish} becomes
\begin{equation*}
H^{2}(X,O_{X})= H^{3}(X,O_{X})=0.
\end{equation*}

If we further assume that $X$ is a surface, then Grothendieck vanishing theorem says that $H^{3}(X,O_{X})=0$. Then the condition \eqref{eq: condi-vanish} becomes
\begin{equation*}
H^{2}(X,O_{X})=0.
\end{equation*}In other words, we require that the geometric genus of $X$ is trivial, which is exactly the hypothesis of Conjecture \ref{c:Bloch}.

In the case of $p=3$, $p-2=1$. Since $0 \leq i \leq p-2$, $i$ can be $0$ or $1$. The condition \eqref{eq: condi-vanish} becomes
\begin{align*}
\begin{cases}
   H^{3}(X,O_{X})= H^{4}(X,O_{X})= H^{5}(X,O_{X})=0,  \\
  H^{3}(X,\Omega^{1}_{X/k})= H^{4}(X,\Omega^{1}_{X/k})=0. 
\end{cases}
\end{align*}
This is the condition used by Mackall to study the pro-representability of $\widehat{CH}^{3}(A)$, cf. Theorem \ref{t:mackall}.

In general, the condition \eqref{eq: condi-vanish} carries the following form:
{\footnotesize
 \begin{align}\label{eq: condi-vanish-explicit}
\begin{cases}
   i=0: H^{p}(X,O_{X}) = H^{p+1}(X,O_{X})= \cdots  = H^{2p-1}(X,O_{X})=0, \\
   i=1: H^{p}(X,\Omega^{1}_{X/k}) = H^{p+1}(X,\Omega^{1}_{X/k})= \cdots =H^{2p-2}(X,\Omega^{1}_{X/k})=0, \\
    \ \ \   \vdots   \ \ \  \ \ \  \ \ \  \ \ \  \ \ \ \vdots  \\
i=p-2: H^{p}(X,\Omega^{p-2}_{X/k}) = H^{p+1}(X,\Omega^{p-2}_{X/k})=0.
\end{cases}
\end{align}
}

It is interesting to find varieties that satisfy the condition \eqref{eq: condi-vanish}. We start with projective spaces, which is classical and well-known.
\begin{example}
Let $X=\mathbb{P}_{k}^{n}$ be the $n$-dimensional projective space. The cohomology $H^{j}(X, \Omega^{j'}_{X/k})$ is computed as
 \begin{align*}
\begin{cases}
   H^{j}(X, \Omega^{j'}_{X/k})=k, \ \mathrm{if} \ j=j'\leq n, \\
H^{j}(X, \Omega^{j'}_{X/k})=0, \ \mathrm{else}.
\end{cases}
\end{align*}Hence, projective spaces obviously satisfy the condition \eqref{eq: condi-vanish}.

\end{example}

 Next we consider suitable hypersurfaces or complete intersections in projective spaces. This guides us to the work by Rapoport.
\begin{example}[\cite{Ra}]
For the case of $p=3$, the following varieties satisfy the condition \eqref{eq: condi-vanish}:
\begin{itemize}
\item (1) a smooth cubic hypersurface in $\mathbb{P}^{6}$.
\item (2) a smooth complete intersection of two quadrics in $\mathbb{P}^{7}$.
\end{itemize}

\end{example}
This example is used in the literature; see, for instance, Example 3.2 of \cite{Ma} for related discussions.

\begin{remark} \label{rem:example}
When the ground field is the complex number field $\mathbb{C}$, part of the condition \eqref{eq: condi-vanish} is that the $p+1$-th cohomology groups are trivial:
\begin{equation} \label{eq:GG}
H^{p+1}(X,O_{X}) = H^{p+1}(X,\Omega^{1}_{X/ \mathbb{C}})= \cdots = H^{p+1}(X,\Omega^{p-2}_{X/ \mathbb{C}})=0.
\end{equation}This is used by Green and Griffiths to study deformations of the Chow group $CH^{p}(X)$ in \cite{GGChow}. Let $Z \subset X$ be a codimension-$p$ subvariety and let $[Z] \in CH^{p}(X)$ be the corresponding rational equivalence class. They prove that, under the condition \eqref{eq:GG} together with assuming that the fundamental class $\psi_{0}(Z)$ remains Hodge type $(p,p)$, the class $[Z]$ can be lifted formally. We refer to \cite{GGChow} for the precise statements and more discussions.

\end{remark}

For any $A \in Art_{k}$, we write $X_{A}=X \otimes_{\mathrm{Spec}(k)} \mathrm{Spec}(A)$ as before. There is an isomorphism, which can be checked locally,
\begin{equation} \label{eq:isoOmega}
\Omega^{l}_{X_{A}/k} \cong \bigoplus_{l_{1}+l_{2}=l}\Omega^{l_{1}}_{X/k}\otimes_{k}\Omega^{l_{2}}_{A/k},
\end{equation}where $l$, $l_{1}$ and $l_{2}$ are non-negative integers. Recall that $i \leq p-2$ in Setting \ref{s:set1}. We let $l=i$ and obtain that
\begin{align} \label{eq:cohzero}
H^{j}(X, \Omega^{i}_{X_{A}/k})& \cong H^{j}(X, \bigoplus_{l_{1}+l_{2}=i}\Omega^{l_{1}}_{X/k}\otimes_{k}\Omega^{l_{2}}_{A/k}) \\ \notag
& \cong \bigoplus_{l_{1}+l_{2}=i}H^{j}(X, \Omega^{l_{1}}_{X/k}\otimes_{k}\Omega^{l_{2}}_{A/k}) \\  \notag
& \cong \bigoplus_{l_{1}+l_{2}=i}H^{j}(X, \Omega^{l_{1}}_{X/k})\otimes_{k}\Omega^{l_{2}}_{A/k}, \notag
\end{align}where $j$ is a non-negative integer.

Note that $0 \leq l_{1} \leq i \leq p-2$. The condition \eqref{eq: condi-vanish} says that
\begin{equation*}
H^{j}(X, \Omega^{l_{1}}_{X/k})=0, \ \mathrm{where} \ j=p, \ p+1, \ \cdots,\ 2p-1-l_{1}.
\end{equation*}Since $l_{1} \leq i$, $2p-1-l_{1} \geq 2p-1-i$, it follows that
\begin{equation} \label{eq:cohzero2}
H^{j}(X, \Omega^{l_{1}}_{X/k})=0, \ \mathrm{where} \ j=p, \ p+1, \ \cdots,\ 2p-1-i.
\end{equation}

As a consequence of \eqref{eq:cohzero} and \eqref{eq:cohzero2}, we obtain the following result.
\begin{lemma} \label{l:vanishext}
In Setting \ref{s:set1}, we have the following
{\footnotesize
\begin{equation*}
H^{p}(X,\Omega^{i}_{X_{A}/ k})= H^{p+1}(X,\Omega^{i}_{X_{A}/ k})=\cdots =H^{2p-1-i}(X,\Omega^{i}_{X_{A}/ k})=0.
\end{equation*}
}

\end{lemma}

The kernel of 
\[
d: O_{X_{A}} \to \Omega^{1}_{X_{A}/k},
\]is the de Rham cohomology sheaf $\mathcal{H}_{dR}^{0}(O_{X_{A}})$, cf. Lemma \ref{l:derhshiso} (let $q=0$). We denote by $dO_{X_{A}}$ the image of $d$. There is a short exact sequence of sheaves (from the fundamental theorem for group homomorphisms)
\begin{equation*}
0 \to \mathcal{H}_{dR}^{0}(O_{X_{A}}) \to O_{X_{A}} \to dO_{X_{A}} \to 0,
\end{equation*}whose associated long exact sequence has the form
\begin{align} \label{eq:lesdox}
& \cdots \to H^{j}(X, \mathcal{H}_{dR}^{0}(O_{X_{A}})) \to H^{j}(X, O_{X_{A}}) \to H^{j}(X, dO_{X_{A}}) \\
& \to H^{j+1}(X, \mathcal{H}_{dR}^{0}(O_{X_{A}})) \to \cdots, \notag
\end{align}where $j$ is a non-negative integer.

When $j > 0$, it follows from Corollary \ref{c:BO-A} that 
\[
H^{j}(X, \mathcal{H}_{dR}^{0}(O_{X_{A}}))=H^{j+1}(X, \mathcal{H}_{dR}^{0}(O_{X_{A}}))=0.
\]The exactness of \eqref{eq:lesdox} yields that 
\[
H^{j}(X, O_{X_{A}}) \cong H^{j}(X, dO_{X_{A}}).
\]

Combining it with Lemma \ref{l:vanishext} (let $i=0$), we obtain the following vanishing result.
\begin{lemma} \label{l:doxA}
In Setting \ref{s:set1}, for $j=p$, $p+1$, $\cdots$, $2p-1$, the $j$-th cohomology group $H^{j}(X, dO_{X_{A}})$ is trivial, i.e.,
\begin{equation*}
H^{p}(X, dO_{X_{A}})=H^{p+1}(X, dO_{X_{A}})=\cdots =H^{2p-1}(X, dO_{X_{A}})=0.
\end{equation*}

\end{lemma}

Let $T^{1}$ and $d\Omega^{1}_{X_{A}/k}$ denote the kernel and the image of the following map respectively
\[
d: \ \Omega^{1}_{X_{A}/k} \to \Omega^{2}_{X_{A}/k}.
\]

There are two short exact sequences of sheaves
\begin{align} \label{eq:SEQ1}
0 \to  dO_{X_{A}} \to T^{1} \to \mathcal{H}_{dR}^{1}(O_{X_{A}}) \to 0, 
\end{align}
\begin{align} \label{eq:SEQ2}
0 \to  T^{1} \to \Omega^{1}_{X_{A}/k} \to d\Omega^{1}_{X_{A}/k} \to 0.
\end{align}Here the exactness of \eqref{eq:SEQ1} is from Lemma \ref{l:derhshiso} (let $q=1$), and the exactness of \eqref{eq:SEQ2} is from the fundamental theorem for group homomorphisms.

The long exact sequence of sheaf cohomology groups associated to \eqref{eq:SEQ1} has the form
{\footnotesize
\begin{align} \label{eq:long-SEQ11}
 \cdots \to H^{j}(X, dO_{X_{A}}) \to H^{j}(X, T^{1}) \to H^{j}(X, \mathcal{H}_{dR}^{1}(O_{X_{A}})) \to \cdots,
\end{align}
}where $j$ is a non-negative integer.

When $j$ is an integer from $p$ to $2p-1$, i.e., $j=p$, $p+1$, $\cdots$, $2p-1$, it has been proved in Lemma \ref{l:doxA} that 
\[
H^{j}(X, dO_{X_{A}})=0.
\]Meanwhile, since $j \geq p \geq 2 >1$, $H^{j}(X, \mathcal{H}_{dR}^{1}(O_{X_{A}}))=0$ by Corollary \ref{c:BO-A}. The exactness of the sequence \eqref{eq:long-SEQ11} shows
\begin{equation} \label{eq:T1zero}
 H^{j}(X, T^{1})=0, \ \mathrm{where} \ j=p,\ p+1,\ \cdots, \ 2p-1.
\end{equation}

The long exact sequence of sheaf cohomology groups associated to \eqref{eq:SEQ2} has the form
\begin{align} \label{eq:long-SEQ22}
 & \cdots \to H^{j}(X, T^{1}) \to H^{j}(X, \Omega^{1}_{X_{A}/k}) \to  H^{j}(X, d\Omega^{1}_{X_{A}/k}) \to \\
& H^{j+1}(X, T^{1}) \to H^{j+1}(X, \Omega^{1}_{X_{A}/k}) \to \cdots,  \notag
\end{align}
where $j$ is a non-negative integer.

For $j=p$, $p+1$, $\cdots$, $2p-2$, substituting \eqref{eq:T1zero} into \eqref{eq:long-SEQ22}, we see that 
\begin{equation*} \label{eq:cohdomega1}
H^{j}(X, \Omega^{1}_{X_{A}/k}) \cong  H^{j}(X, d\Omega^{1}_{X_{A}/k}).
\end{equation*}Combining it with Lemma \ref{l:vanishext} (let $i=1$), we obtain the following vanishing result.
\begin{lemma} \label{l:domega1}
For $j=p$, $p+1$, $\cdots$, $2p-2$, the cohomology $H^{j}(X, d\Omega^{1}_{X_{A}/k})$ is trivial. That is, 
\[
H^{p}(X, d\Omega^{1}_{X_{A}/k})=H^{p+1}(X, d\Omega^{1}_{X_{A}/k})= \cdots =H^{2p-2}(X, d\Omega^{1}_{X_{A}/k})=0.
\]

\end{lemma}

We now generalize Lemma \ref{l:doxA} and Lemma \ref{l:domega1}. Let $T^{i}$ and $d\Omega^{i}_{X_{A}/k}$ denote the kernel and the image of the following map respectively
\[
d: \ \Omega^{i}_{X_{A}/k} \to \Omega^{i+1}_{X_{A}/k},
\]where $i=0$, $1$, $\cdots$, $p-2$. 

\begin{lemma} \label{l:domega-n}
In Setting \ref{s:set1}, for $j=p$, $p+1$, $\cdots$, $2p-1-i$, the cohomology $H^{j}(X, d\Omega^{i}_{X_{A}/k})$ is trivial. In other words, we have
\[
H^{p}(X, d\Omega^{i}_{X_{A}/k})=H^{p+1}(X, d\Omega^{i}_{X_{A}/k})= \cdots =H^{2p-1-i}(X, d\Omega^{i}_{X_{A}/k})=0.
\]

In particular, when $i=p-2$, we have
\begin{equation} \label{eq:use}
H^{p}(X, d\Omega^{p-2}_{X_{A}/k})=H^{p+1}(X, d\Omega^{p-2}_{X_{A}/k})=0.
\end{equation}
\end{lemma}

\begin{proof}
We use induction on $i$. The case of $i=0$ has been proved in Lemma \ref{l:doxA}. For the case of $i \geq 1$, suppose that the case of $i-1$ has been proved. That is, for 
\[
j=p, \ p+1, \  \cdots, \ 2p-1-(i-1),
\]the cohomology $H^{j}(X, d\Omega^{i-1}_{X_{A}/k})$ is trivial, i.e.,
\begin{equation} \label{eq:induc}
H^{j}(X, d\Omega^{i-1}_{X_{A}/k})=0.
\end{equation}

We use the method for proving Lemma \ref{l:domega1} to proceed. There are two short exact sequences of sheaves
\begin{align} \label{eq:SEQ1i}
0 \to  d\Omega^{i-1}_{X_{A}/k} \to T^{i} \to \mathcal{H}_{dR}^{i}(O_{X_{A}}) \to 0, 
\end{align}
\begin{align} \label{eq:SEQ2i}
0 \to  T^{i} \to \Omega^{i}_{X_{A}/k} \to d\Omega^{i}_{X_{A}/k} \to 0.
\end{align}Here the exactness of \eqref{eq:SEQ1i} is from Lemma \ref{l:derhshiso} (let $q=i$), and the exactness of \eqref{eq:SEQ2i} is from the fundamental theorem for group homomorphisms.

The long exact sequence of sheaf cohomology groups associated to \eqref{eq:SEQ1i} has the form
{\footnotesize
\begin{align} \label{eq:long-SEQ1i}
\cdots \to  H^{j}(X, d\Omega^{i-1}_{X_{A}/k}) \to H^{j}(X, T^{i}) \to H^{j}(X, \mathcal{H}_{dR}^{i}(O_{X_{A}})) \to \cdots. 
\end{align}
}

Since the integer $i$ satisfies that $0 \leq i \leq p-2$, for $j=p, \ p+1, \  \cdots, \ 2p-1-(i-1)$, it is clear that $j \geq p > p-2 \geq i$. It follows from Corollary \ref{c:BO-A} that 
\[
H^{j}(X, \mathcal{H}_{dR}^{i}(O_{X_{A}}))=0.
\]Combining with \eqref{eq:induc}, we use the exactness of the sequence \eqref{eq:long-SEQ1i} to see that 
\begin{equation} \label{eq:Ti}
H^{j}(X, T^{i})=0, \ \mathrm{where} \ j=p, \ p+1, \ \cdots, \ 2p-1-(i-1).
\end{equation}

The long exact sequence of sheaf cohomology groups associated to \eqref{eq:SEQ2i} has the form
\begin{align} \label{eq:long-SEQ2}
& \cdots \to  H^{j}(X, T^{i}) \to H^{j}(X, \Omega^{i}_{X_{A}/k}) \to H^{j}(X, d\Omega^{i}_{X_{A}/k})  \\
& \to H^{j+1}(X, T^{i}) \to H^{j+1}(X, \Omega^{i}_{X_{A}/k}) \to \cdots. \notag
\end{align}

When $j=p$, $p+1$, $\cdots$, $2p-1-i$, it follows from \eqref{eq:Ti} that, both $ H^{j}(X, T^{i})$ and $ H^{j+1}(X, T^{i})$ are trivial. The exactness of \eqref{eq:long-SEQ2} yields that
\[
 H^{j}(X, \Omega^{i}_{X_{A}/k}) \cong H^{j}(X, d\Omega^{i}_{X_{A}/k}).
\]Combining with Lemma \ref{l:vanishext}, we obtain that 
\[
 H^{j}(X, d\Omega^{i}_{X_{A}/k})=0, \ \mathrm{where} \ j=p, \ p+1, \ \cdots, \ 2p-1-i.
\]This completes the proof.

\end{proof}

There is a short exact sequence of sheaves
\begin{align*} 
0 \to  d\Omega^{p-2}_{X_{A}/k} \to \Omega^{p-1}_{X_{A}/k} \to \frac{\Omega^{p-1}_{X_{A}/k}}{d\Omega^{p-2}_{X_{A}/k}} \to 0,
\end{align*}whose associated long exact sequence has the form
\begin{align*} 
& \cdots \to H^{p}(X,d\Omega^{p-2}_{X_{A}/k}) \to H^{p}(X,\Omega^{p-1}_{X_{A}/k}) \to H^{p}(X,\frac{\Omega^{p-1}_{X_{A}/k}}{d\Omega^{p-2}_{X_{A}/k}}) \\
& \to  H^{p+1}(X,d\Omega^{p-2}_{X_{A}/k}) \to \cdots. \notag
\end{align*}Since $H^{p}(X,d\Omega^{p-2}_{X_{A}/k})= H^{p+1}(X,d\Omega^{p-2}_{X_{A}/k})=0$, cf. \eqref{eq:use}, it follows that
\begin{align} \label{eq:iso}
H^{p}(X,\frac{\Omega^{p-1}_{X_{A}/k}}{d\Omega^{p-2}_{X_{A}/k}}) & \cong H^{p}(X,\Omega^{p-1}_{X_{A}/k}) \\ \notag
& \cong H^{p}(X,\bigoplus_{l_{1}+l_{2}=p-1}\Omega^{l_{1}}_{X/k}\otimes_{k} \Omega^{l_{2}}_{A/k}) \\ \notag
& \cong \bigoplus_{l_{1}+l_{2}=p-1} H^{p}(X,\Omega^{l_{1}}_{X/k}\otimes_{k} \Omega^{l_{2}}_{A/k}) \\ \notag
& \cong \bigoplus_{l_{1}+l_{2}=p-1} H^{p}(X,\Omega^{l_{1}}_{X/k}) \otimes_{k} \Omega^{l_{2}}_{A/k},   \notag
\end{align}where the second isomorphism is from  \eqref{eq:isoOmega} (let $l=p-1$).

When $l_{1}$ is from $0$ to $p-2$, i.e., $0 \leq l_{1} \leq p-2$, the condition \eqref{eq: condi-vanish} says that 
\[
H^{p}(X,\Omega^{l_{1}}_{X/k})=0.
\]This shows that the only non-trivial direct summand of the above direct sum is $H^{p}(X,\Omega^{p-1}_{X/k}) \otimes_{k} A$, in which case $l_{1}=p-1$ and $l_{2}=0$. As a consequence, the isomorphism \eqref{eq:iso} becomes 
\begin{equation} \label{eq:iso2}
H^{p}(X,\frac{\Omega^{p-1}_{X_{A}/k}}{d\Omega^{p-2}_{X_{A}/k}}) \cong H^{p}(X,\Omega^{p-1}_{X/k}) \otimes_{k} A.
\end{equation}

When $A=k$, the isomorphism \eqref{eq:iso2} specializes to 
\begin{equation} \label{eq:iso3}
H^{p}(X,\frac{\Omega^{p-1}_{X/k}}{d\Omega^{p-2}_{X/k}}) \cong H^{p}(X,\Omega^{p-1}_{X/k}) \otimes_{k} k \cong H^{p}(X,\Omega^{p-1}_{X/k}).
\end{equation}

By Lemma \ref{lemma: l-l-omega}, we obtain the following isomorphisms of sheaves \footnote{These cyclic homologies here are defined over $\mathbb{Q}$.}
\begin{align}\label{eq: HC-smooth-list}
 HC^{(p-1)}_{p-1}(O_{X_{A}}) \cong \frac{\Omega^{p-1}_{X_{A}/\mathbb{Q}}}{d\Omega^{p-2}_{X_{A}/\mathbb{Q}}}, HC^{(p-1)}_{p-1}(O_{X}) \cong \frac{\Omega^{p-1}_{X/\mathbb{Q}}}{d\Omega^{p-2}_{X/\mathbb{Q}}}.
\end{align}

From now on, we put an additional assumption on the ground field $k$ and consider the following setting.
\begin{setting} \label{s:set2}
In Setting \ref{s:set1}, we further assume that $k$ is an algebraic field extension of $\mathbb{Q}$.

\end{setting}

In this setting, for each non-negative integer $l$, $\Omega^{l}_{X_{A}/\mathbb{Q}}=\Omega^{l}_{X_{A}/k}$. Hence, we can rewrite \eqref{eq: HC-smooth-list} as
\begin{align}\label{eq: HC-smooth-list2}
 HC^{(p-1)}_{p-1}(O_{X_{A}}) \cong \frac{\Omega^{p-1}_{X_{A}/k}}{d\Omega^{p-2}_{X_{A}/k}}, HC^{(p-1)}_{p-1}(O_{X}) \cong \frac{\Omega^{p-1}_{X/k}}{d\Omega^{p-2}_{X/k}}.
\end{align}

Let $\overline{HC}^{(p-1)}_{p-1}(O_{X_{A}})$ be the relative cyclic homology sheaf, cf. Definition \ref{d:relcyc} (let $i=n=p-1$). There is a split short exact sequence
\[
0 \to \overline{HC}^{(p-1)}_{p-1}(O_{X_{A}}) \to HC^{(p-1)}_{p-1}(O_{X_{A}}) \xrightarrow{A \to k} HC^{(p-1)}_{p-1}(O_{X}) \to 0,
\]which yields a split short exact sequence (by \eqref{eq: HC-smooth-list2}):
\begin{equation*}
0 \to \overline{HC}^{(p-1)}_{p-1}(O_{X_{A}}) \to \frac{\Omega^{p-1}_{X_{A}/k}}{d\Omega^{p-2}_{X_{A}/k}} \xrightarrow{A \to k} \frac{\Omega^{p-1}_{X/k}}{d\Omega^{p-2}_{X/k}} \to 0.
\end{equation*}

Passing to cohomology groups, we have a short exact sequence
{\footnotesize
\[
 0 \to H^{p}(X,\overline{HC}^{(p-1)}_{p-1}(O_{X_{A}})) \to H^{p}(X,\frac{\Omega^{p-1}_{X_{A}/k}}{d\Omega^{p-2}_{X_{A}/k}}) \xrightarrow{A \to k} H^{p}(X,\frac{\Omega^{p-1}_{X/k}}{d\Omega^{p-2}_{X/k}}) \to 0,
\]
}which can be rewritten as (by using \eqref{eq:iso2} and \eqref{eq:iso3})
{\footnotesize
\[
 0 \to H^{p}(X,\overline{HC}^{(p-1)}_{p-1}(O_{X_{A}})) \to H^{p}(X,\Omega^{p-1}_{X/k}) \otimes_{k} A \xrightarrow{A \to k} H^{p}(X,\Omega^{p-1}_{X/k}) \to 0.
\]
}
Since $A=k \oplus m_{A}$ (as $k$-vector spaces), it follows that
\begin{align*} 
H^{p}(X,\Omega^{p-1}_{X/k}) \otimes_{k} A  & \cong  H^{p}(X,\Omega^{p-1}_{X/k}) \otimes_{k}(k \oplus m_{A}) \\ \notag
& \cong  (H^{p}(X,\Omega^{p-1}_{X/k}) \otimes_{k} k) \oplus (H^{p}(X,\Omega^{p-1}_{X/k}) \otimes_{k} m_{A})\\ \notag
& \cong  H^{p}(X,\Omega^{p-1}_{X/k}) \oplus (H^{p}(X,\Omega^{p-1}_{X/k}) \otimes_{k} m_{A}). \notag
\end{align*}

As a consequence, we have
\[
H^{p}(X,\overline{HC}^{(p-1)}_{p-1}(O_{X_{A}})) \cong H^{p}(X,\Omega^{p-1}_{X/k}) \otimes_{k} m_{A}.
\]
Combining it with Lemma \ref{l:chhatisocyc}, we obtain the main result of this paper.
\begin{theorem} \label{t:yang}
In Setting \ref{s:set2}, for any $A \in Art_{k}$, there is an isomorphism
\[
\widehat{CH}^{p}(A) \cong H^{p}(X, \Omega^{p-1}_{X/ k})\otimes_{k}m_{A}.
\]Consequently, the functor $\widehat{CH}^{p}$ is pro-representable, cf. \eqref{eq: pro-rep vs}.

\end{theorem}

\textbf{Acknowledgments}. The author is very grateful to Spencer Bloch \cite{Bl5} and Jan Stienstra \cite{Stien2} for sharing their ideas. He thanks Aise Johan de Jong, Phillip Griffiths, Jerome William Hoffman, Kefeng Liu, Yang Shen, Wei Xia, Chenglong Yu and ZhiweiZheng for discussions and /or comments.

This work is partially supported by the Key Program and the Major Program of Natural Science Research Foundation of Anhui Provincial Education Department (No.2024AH051403 and No.2023AH040225), the Research and Innovation Team of Anhui Province (No.2024AH010023), and the Start-up Research Fund of Chuzhou University (No.2024qd53).


\end{document}